\input amstex
\documentstyle{amsppt}
\document

\magnification=1200

\topmatter
\title  {A note on a theorem of Bowditch }
\endtitle
\rightheadtext{A note on a theorem of Bowditch }

\author{ Panos Papasoglu }
\endauthor

\date { }
\enddate

\address{
Departement de Math\'ematiques, Universit\'e de Paris XI (Paris-Sud)
91405 Orsay, FRANCE.}
\endaddress
\email{  panos\@ math.u-psud.fr}
\endemail

\endsubjclass
\abstract { Bowditch showed that a one-ended hyperbolic group which is not
a triangle group
splits over
a two-ended group if and only if its boundary has a local cut point. As a
corollary one obtains that splittings of hyperbolic groups over two-ended
groups are preserved under quasi-isometries. In this note we give a more direct
proof of this corollary.}

\endabstract

\endtopmatter

\NoBlackBoxes
\heading {\bf \S 0. Introduction } \endheading
In this paper we give a new proof Bowditch's result (\cite {Bo1}) that splittings of one-ended
hyperbolic groups over two ended groups are preserved by quasi-isometries.
Our proof is similar to the proof of Stallings' theorem on groups with infinitely
many ends (see for example \cite {S-W}). The similarity of these two theorems
suggests several questions that we present at the end of this paper.\newline
We recall that Bowditch showed (\cite {Bo2,Bo3}) that if the boundary, $\partial G$, of a
one-ended hyperbolic group $G$ has a cut point then $G$ splits over a 2-ended group. Swarup
(\cite {Sw }) using Bowditch's methods subsequently showed that $\partial G$ has no 
cut points. This implies  that $\partial G$ is locally connected by earlier results
of Bestvina-Mess (\cite {B-M}).\newline
The proof that a one-ended hyperbolic group whose boundary has a cut point splits
over a 2-ended group, has 3 steps:\newline
1. Pass from the set of cut points of $\partial G$ to an $\Bbb R$-tree $T$ on which
$G$ acts by homeomorphisms. There is a nice account of this step in \cite {Swe }. We recall
that Whyburn (\cite {Wh }) was the first to study the cut point set of a locally
connected metric space and to note that it is `treelike' (a dendrite).\newline
2. Note that one can change the action of $G$ on $T$ to an action by isometries
and two-ended segment stabilizers. An elegant way to do this is given in \cite {Le}.\newline
3. Apply Rips' machine (\cite {B-F}) to conclude that $G$ splits over a two-ended group.\newline
We note that a recent result of Delzant-Potyagailo (\cite {D-P}) can also be used to
infer from this splitting theorem that in fact the boundary of a 1-ended hyperbolic group
has no cut points. In any case to establish quasi-isometry invariance of splittings
of hyperbolic groups over 2-ended groups one does not need this stronger result.\newline
In this paper we show the following:
\proclaim {Theorem}
Let $G$ be a one-ended hyperbolic group. Suppose that a pair of points $\{x,y \}$
separates $\partial G$ (i.e. $\partial G- \{x,y \}$ is not connected). Then either
$\partial G$ is a circle or $G$ splits over a 2-ended group.
\endproclaim
One easily deduces that splittings over 2-ended groups are invariant under quasi-isometries:
\proclaim {Corollary}
Let $G,H$ be  one-ended quasi-isometric hyperbolic groups. Suppose that $G$
is not virtually a surface group and that it splits over a 2-ended group.
Then $H$ splits over a 2-ended group.
\endproclaim
\demo {Proof }
We note that quasi-isometric groups have homeomorphic boundaries. If $\partial G$
has a cut point then $\partial H$ has a cut point so $H$ splits over a 2-ended group.
We suppose therefore that $\partial G, \partial H$ have no cut points. Since $G$
splits over a a 2-ended group $D$ , $\partial D$ separates $\partial G$ (to see this
note that a neighborhood of $D$ separates the Cayley graph of $G$). Since
$D$ is 2-ended $\partial D$ has two points. By the result of Tukia-Gabai
 (\cite {T,Ga,C-J}) $\partial G$ is a circle if and only if $G$ is virtually
a surface group. We note here that surface groups split over $\Bbb Z$ but
hyperbolic triangle groups which are commensurable (and hence quasi-isometric) to
surface groups don't split over 2-ended groups, and this is why we exclude such
groups in the statement of the theorem. Since $\partial H$ is not a circle
,has no cut points, and is separated by 2 points, we can apply the theorem above
to conclude that $H$ splits over a 2-ended group.
\enddemo
We note that in this paper we don't recover all results of \cite {Bo1}. In particular
in \cite {Bo1} is shown that $G$ splits if $\partial G$ has a {\it local} cut point. This
is a deeper result (see \cite {Bo1}, sec. 5).\newline
We will use the fact that the action of a hyperbolic group $G$ on its boundary, $\partial G$,
is a {\it convergence group action}. We recall here what this means:\newline
Suppose that $(g_n)$ is an infinte sequence of distinct elements of $G$. Then
there is a subsequence $(g_i)$
and $x,y\in \partial G$ such that the maps $g_i|(\partial G-\{y \})$
converge uniformly on compact subsets of $\partial G-\{y \}$ to $x$.
\heading {\bf \S 1. Quasi-isometry invariance of splittings } \endheading
We will need the following simple lemma:
\proclaim {Lemma 1}
Let $X$ be a compact, path connected, metric space with no cut points.
Suppose that 
for any pair of points $\{a,b\}$ of $X$ there is a pair of points $\{x,y \}$ such that
$a,b$ lie in distinct components of $X-\{x,y \}$.
Then $X$ is homeomorphic to $S^1$.
\endproclaim
\demo {Proof}
Since $X$ has no cut points $X$ contains a simple closed curve $p$. We will show
that $X=p$. Indeed suppose that $x\in X-p$. Let $q$ be a simple path joining $x$ to $p$
and let $y=q\cap p$. 
Since $X$ has no cut points there is
a simple path $q_1$ joining $x$ to $p$ such that $y\notin q_1$. Let $y_1=p\cap q_1$.
It is clear now that for any pair of points $x_1,x_2\in X$, $y,y_1$ lie in the same
component of $X-\{x_1,x_2\}$.
\enddemo
\proclaim {Theorem}
Let $G$ be a one-ended hyperbolic group. Suppose that a pair of points $\{x,y \}$
separates $\partial G$ (i.e. $\partial G- \{x,y \}$ is not connected). Then either
$\partial G$ is a circle or $G$ splits over a 2-ended group.
\endproclaim
\demo {Proof}
By lemma 1 if $\partial G$ is not a circle there are $a,b\in \partial G$ such that
$\{a,b\}$ can not be separated by any pair $t,s\in \partial G$. We pose $X=\partial G$.
We will say that a closed set $C$ is a {\it slice } of $X$ if
 $C$
is the closure of a connected component of $X-\{ s,t \}$ for some $s,t\in X$.
Clearly any slice of $X$ containing $a$ contains also $b$. We consider
all slices of $X$ containing $a$. We have the following lemma:
\proclaim {Lemma 2}
Let $C_1\supset C_2\supset...$ be a descending sequence of slices of $X$.
If $C=\cap C_i$ is non empty then $C$ is a slice of $X$.
\endproclaim
\demo {Proof }
Let's say that $C_i$ is the closure of a component of $X-\{s_i,t_i\} $.
By passing to a subsequence we can assume that $s_i\to s$ , $t_i\to t$ and
$t\ne s$.
$C$ is then the closure of a connected component of $X-\{ s,t\} $ therefore
it is a slice of $X$.
\enddemo
By the lemma above it follows that there is a minimal slice of $X$ containing
$a$. Suppose then that $C$ is a minimal slice of $X$ containing $a$.
Let's say that $C$ is the closure of a component of $X-\{ s,t\} $. We have the following lemma:
\proclaim {Lemma 3}
For any $g\in G$, $gs,gt$ lie in the same component of $X-\{s,t\} $.
\endproclaim
\demo {Proof}
Suppose that the lemma is false. Then $s,t$ lie in different components of
$X-\{g^{-1}s,g^{-1}t\} $. This however contradicts the minimality of $C$.
\enddemo
We consider now the set $S=\{(gs,gt):g\in G\} $. We say that
$(s_1,t_1)\in S$ separates $(s_2,t_2),(s_3,t_3)\in S$ if $s_2,t_2,s_3,t_3$ are not all 
contained in the same slice of $X-\{ s_1,t_1 \} $. We call $(s_1,t_1),(s_2,t_2)\in S$
adjacent if no $(s_3,t_3)$ in $S$ separates them.
We have the following:
\proclaim {Lemma 4}
If $(s_1,t_1),(s_2,t_2)\in S$ then there are finitely many couples $(s',t')\in S$
such that $(s',t') $ separates $(s_1,t_1),(s_2,t_2)$.
\endproclaim
\demo {Proof}
Suppose there is an infinite sequence $(g_ns,g_nt)\in S $ separating 
$(s_1,t_1),(s_2,t_2)$. By passing to a subsequence we can assume
that $g_ns\to s'$, $g_nt\to t'$. Since $X$ has no cut points $s'\ne t'$.
By passing to a subsequence we have that there are points $s",t"$ such
that $g_n\, |X-\{t" \}$ converges on compact sets to $s"$. It follows
that either $s'=s"$ or $t'=s"$. Both of these are impossible since
$g_nt\to t'$ and $g_ns\to s'$.
\enddemo
We would like to join adjacent elements of $S$  by edges and claim that the graph
obtained is a tree. This is not quite true since it is possible
that 3 elements of $S$ are adjacent to each other. Think of a tripod, if one takes 
as $S$ only the endpoints of the tripod and defines separation in the obvious way
it is clear that in order to obtain a tree from $S$ one has to add a new vertex (the
middle vertex of the tripod). We 'complete' then $S$ to a vertex set of a tree as
follows:\newline
We call a subset of $S$, $V$, a vertex set if any two $x,y\in V$ are adjacent
and $V$ is maximal with this property. We define now a graph $T$ as follows:
The vertices of $T$ are the elements of $S$ and the vertex sets of $S$. We join
two vertices of $T$ by an edge if one of them is an element of $S$ and the other
is a vertex set containing it. Note that when $S$ is the vertex set of a tree the
graph obtained by this operation is the barycentric subdivision of this tree.\newline
Clearly $T$ is a tree and $G$ acts on $T$. If $g\in G$ stabilizes an edge of $T$
then $g$ fixes the endpoint of this edge that lies in $S$ so $g$ fixes a pair of points
$(s_1,t_1)\in S$. This implies that edge stabilizers are 2-ended (recall that $G$ is one-ended).
By the definition of $T$ is clear that there is no global fixed point for this action.
This concludes the proof of the theorem. 

\enddemo
\heading {\bf \S 2. Questions } \endheading
One wonders if one can generalize Bowditch's results to splittings over more
complicated subgroups. Here are some specific questions:\newline
1. Let $G$ be a hyperbolic group such that $\partial G$ is not separated
by any Cantor set. Is it true that $\partial G$ is not separated by a segment?\newline
2. Let $G$ be a one ended hyperbolic group. Suppose that $dim (\partial G)\geq 2$
and that a Cantor set separates $\partial G$. Is it true that $G$ splits over
a virtually free or over a virtual surface group? More generally if $dim (\partial G)\geq n$
and $\partial G$ is separated by a closed subset of dimension $\leq n-2$ is it true
that $G$ splits? (this question has also been asked by Gromov in \cite {Gr }).\newline
We note that Bowditch's result has been extended to all finitely presented groups in
\cite {P}. One can ask the questions above in the more general setting of finitely
presented groups by replacing dimension and separation in the boundary by asymptotic
dimension and separation in the Cayley graph by uniformly embedded subsets.

\enddocument

\bye